\documentclass[11pt,reqno]{amsart}

\usepackage{amsmath,amssymb,amsthm, amsfonts}   
\usepackage{bbm} 
\usepackage{graphicx}   
\usepackage{pstricks}   
\usepackage[colorlinks=true, linkcolor=blue, citecolor=magenta]{hyperref} 
\usepackage{tasks}
\usepackage{enumerate}
\usepackage{caption,subcaption} 

\usepackage{mathtools}      
 \mathtoolsset{showonlyrefs} 

\usepackage{mathrsfs}
\usepackage{multicol}   
\usepackage{fullpage}   

\newtheorem{thm}{Theorem}[section]

\newtheorem{prop}[thm]{Proposition}
\newtheorem{lem}[thm]{Lemma}

\newtheorem*{acknowledgements*}{Acknowledgements}

\theoremstyle{definition}
\newtheorem{defin}[thm]{Definition}

\theoremstyle{remark}

\title{Some geometric properties of Riemann's non-differentiable function}
\author{Daniel Eceizabarrena}
\address{BCAM - Basque Center for Applied Mathematics. Alameda de Mazarredo 14, 48009 Bilbao, Spain.  e-mail: {\tt deceizabarrena@bcamath.org}}
\date{\today}

\begin{document}

\begin{abstract}
Riemann's non-differentiable function is a celebrated example of a continuous but almost nowhere differentiable function. There is strong numeric evidence that one of its complex versions represents a geometric trajectory in experiments related to the binormal flow or the vortex filament equation. In this setting, we analyse certain geometric properties of its image in $\mathbb{C}$. 
The objective of this note is to assert that the Hausdorff dimension of its image is no larger than 4/3 and that it has nowhere a tangent. 
\end{abstract}

\maketitle


\section{Introduction}\label{S_Intro}
Riemann's non-differentiable function
\begin{equation}\label{RiemannFunction}
R(x) = \sum_{n=1}^{\infty}{\frac{ \sin{ (n^2x ) } }{ n^2 }}, \qquad x \in \mathbb{R},
\end{equation}
is a classic example of a continuous but almost nowhere differentiable function. It was proposed by Riemann in the 1860s and collected by Weierstrass in his famous speech \cite{Weierstrass} in the Prussian Academy of Sciences in Berlin in 1872. Since then, it has been studied analytically and in various forms \cite{Duistermaat,Gerver,Gerver2,Hardy,HardyLittlewood,Jaffard} and some results concerning its graph have been obtained  \cite{ChamizoCordoba,HolschneiderTchamitchian}. Generalisations in the exponents of the phase and the denominator of $R$ have also been analysed \cite{ChamizoUbis2007,ChamizoUbis2014,Pastor}. Surprisingly, recent studies suggest that $R$ has a intrinsic geometric nature related to physical experiments. De la Hoz and Vega \cite{DeLaHozVega} showed that a complex generalisation, 
\begin{equation}\label{Phi}
\phi(t) = \sum_{k \in \mathbb{Z}}^{\infty}{ \frac{ e^{-4\pi^2 i k^2 t} - 1 }{-4\pi^2 k^2} }, \qquad t \in \mathbb{R},
\end{equation}
is related to the temporal trajectory of a vortex filament in the context of the binormal flow. This flow is modelled by the Vortex Filament Equation, an evolution equation for a curve $\boldsymbol{X}:\mathbb{R}^2 \to \mathbb{R}^3$, $\boldsymbol{X} = \boldsymbol{X}(s,t)$,  given by 
\begin{equation}\label{VFE}
\boldsymbol{X}_t = \boldsymbol{X}_s \times \boldsymbol{X}_{ss}, \qquad \boldsymbol{X}(s,0) = \boldsymbol{X}_0(s),
\end{equation}
where $t$ represents time and $s$ is the arclength parameter. Assume that $\boldsymbol{X}_M$ is the solution corresponding to the initial datum given by a planar regular polygon of $M \in \mathbb{N}$ sides. Then, strong numeric evidence was given that $\boldsymbol{X}_M(0,t)$ converges to $\phi$ when $M \to \infty$, after some proper rescaling. Hence, the study of the geometry of Riemann's non-differentiable function has a physical justification and goes beyond the mere mathematical curiosity.

\begin{figure}[h]
\centering
\begin{subfigure}{0.475\textwidth}
  \centering
  \includegraphics[width=1\linewidth]{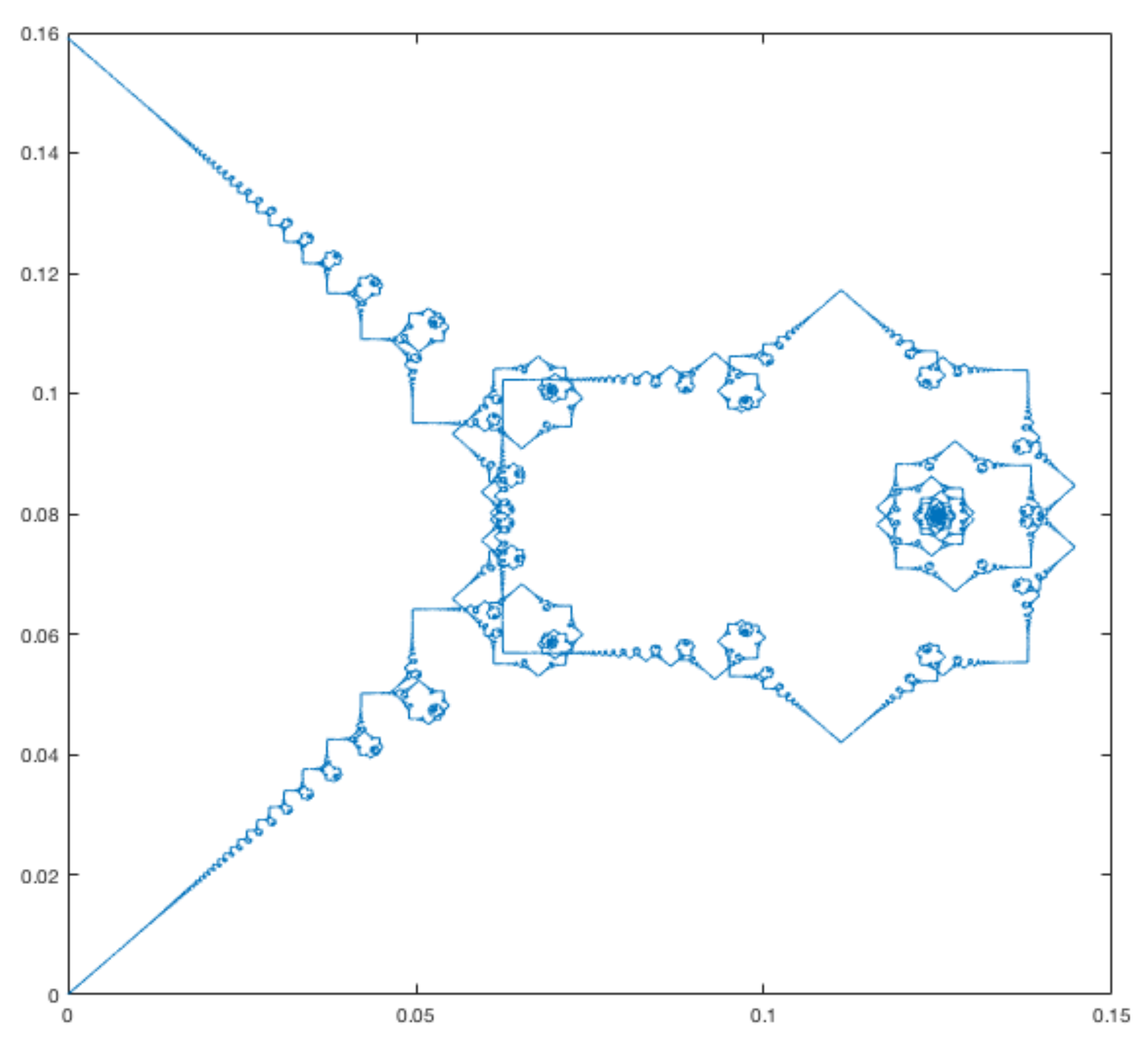}
  \caption{$\phi(t)$}
  \label{FigurePhi}
\end{subfigure}%
\begin{subfigure}{0.05\textwidth}
  \centering
  \hspace{1pt}
\end{subfigure}%
\begin{subfigure}{0.465\textwidth}
  \centering
  \includegraphics[width=1\linewidth]{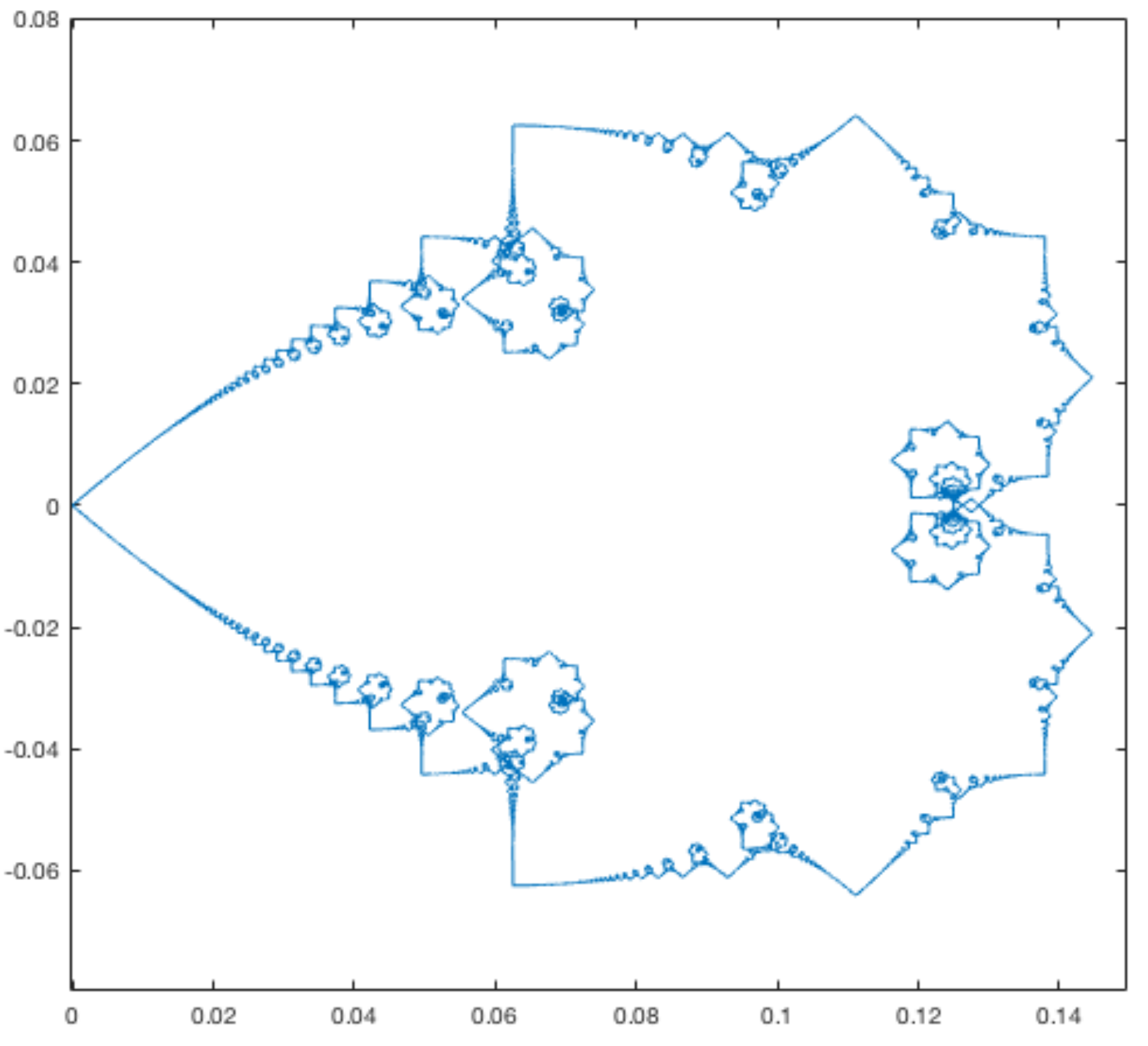}
  \caption{$\phi(t) - it = - i \,\phi_D(-4\pi t)/2\pi + 1/12$}
  \label{FigurePhiDuistermaat}
\end{subfigure}
\caption{Comparison between the images of $[0,1/(2\pi)]$ by $\phi$ and $\phi_D$ as subsets of $\mathbb{C}$.
}
\label{FigurePhiAndPhiDuistermaat}
\end{figure}

Duistermaat \cite{Duistermaat} studied a simpler complex generalisation of $R$,
\begin{equation}\label{PhiDuistermaat}
\phi_D(t) = \sum_{n=1}^{\infty}{ \frac{e^{i\pi n^2 t}}{i\pi n^2} }, \qquad \qquad \operatorname{Re}{\phi_D(t)} =  \frac{1}{\pi} R(\pi t) ,
\end{equation}  
which is 2-periodic and satisfies $\phi(t) = - i \,\phi_D(-4\pi t)/2\pi + it + 1/12$. He gave its asymptotic behaviour around points corresponding to rational numbers. While it is clear that $\phi$ and $\phi_D$ share most analytic properties, the binormal flow experiment shows that $\phi$ and its image in the complex plane are the natural objects to study geometrically instead of $\phi_D$. The difference between the image of both functions can be appreciated in Figure~\ref{FigurePhiAndPhiDuistermaat}. 

In this note, we assert that the image $\phi(\mathbb{R})$ has a Hausdorff dimension not greater than 4/3 and that when interpreted as a trajectory, it has nowhere a tangent. The reader may check \cite{Eceizabarrena1} and \cite{Eceizabarrena2} for the detailed proofs. As a preliminary comment, we remark that unlike $\phi_D$, $\phi$ is not periodic but satisfies $\phi(t + 1/(2\pi)) = \phi(t) + i/(2\pi)$. Thus, it is enough to analyse the set $\phi([0,1/(2\pi)])$, and it is convenient to scale the variable as $t = t_x = x/(2\pi)$ where $x \in (0,1)$.

\section{Hausdorff dimension}\label{S_HausdorffDimension}

As one can see in Figure~\ref{FigurePhi}, the set $\phi([0,1/(2\pi)])$ has self-similar features, which is usually an indicator of fractality. A non-trivial upper bound for its Hausdorff dimension is given in the following theorem, whose complete proof can be checked in \cite{Eceizabarrena1}.

\begin{thm}\label{TheoremDimension}
Let $\phi$ be Riemann's non-differentiable function \eqref{Phi}. Then, 
\[ 1 \leq \operatorname{dim}_{\mathcal{H}}\phi(\mathbb{R}) \leq \frac43. \]
\end{thm}

While the lower bound follows because $\phi$ is a continuous and non-constant curve, the proof of the upper bound consists in covering the set of image points $\phi(t_{\rho})$ corresponding to irrationals $\rho \in \mathbb{R}\setminus \mathbb{Q}$. If $\{ p_n/q_n \}$ are the continued fraction approximations to $\rho$, then it is well-known 
that $|\rho - p_n/q_n| < q_n^{-2}$. This error matches the range of the following estimate, a consequence of the asymptotic behaviour of $\phi$ around rationals which can be essentially deduced from \cite{Duistermaat}.
\begin{lem}\label{LemmaBoundForContinuedFractions}
There exists a constant $C >0$ such that for every irreducible fraction $p/q$, 
\[ \left| \phi(t_{p/q} + h) - \phi(t_{p/q}) \right| \leq C\,\frac{\sqrt{|h|}}{\sqrt{q}} , \qquad \qquad \text{ whenever } \, |h| \leq \frac{1}{q^2}.\]
\end{lem}
Choose then $h = \rho - p_n/q_n$ so that the image of the set $\mathcal{I} = \frac{1}{2\pi}\left( (0,1) \cap \mathbb{I} \right)$ is covered by
\begin{equation}\label{CoverForIrrationals}
\phi(\mathcal{I}) \subset  \bigcup_{\substack{1 \leq p < q \\ \operatorname{gcd}(p,q)=1\\ q \geq Q_0}}B\left( \phi(t_{p/q}),\frac{C}{q^{3/2}} \right), \qquad \forall Q_0 \in \mathbb{N}.
\end{equation}
Let $\alpha >0$ and $\delta >0$. Defining $Q_\delta = \lceil (C/\delta)^{2/3}\rceil + 1$, the above cover gives an upper bound
\[ \mathcal{H}^{\alpha}_{\delta}(\phi(\mathcal{I})) \leq C^{\alpha}\, \sum_{q=Q_{\delta}}^{\infty}{\frac{1}{q^{3\alpha/2-1}}}  \]
for the auxiliary quantities $\mathcal{H}^{\alpha}_{\delta}$. If $\alpha > 4/3$, the geometric series converges, so taking limits $\delta \to 0$ implies that the $\alpha$-Hausdorff measure is $\mathcal{H}^{\alpha}(\phi(\mathcal{I})) = 0$. Consequently, the dimension can be no larger than 4/3.

The same procedure leads to a more general result connected with the remarkable results of Jaffard \cite{Jaffard}. Recall we say $\phi$ is $\alpha$-H\"older in a point $x_0$, and write $\phi \in C^{\alpha}(x_0)$, when there is a polynomial $P$ of degree at most $\alpha$ such that $ \left| \phi(x_0 + h) - P(h) \right| \leq C\,|h|^{\alpha}$ for small enough $h$ and some $C>0$. Define $\alpha(x_0)$ the H\"older exponent of $\phi$ at $x_0$ and the sets $D_{\beta}$ where $\phi$ has exponent $\beta$ as
\[   \alpha(x_0) = \sup\{ \alpha \mid \phi \in C^{\alpha}(x_0) \}, \qquad \qquad  D_{\beta} = \{ x \in \mathbb{R} \mid \alpha(x) = \beta \} .  \]
Jaffard \cite{Jaffard} proved that $\phi$ is a multifractal function in the sense that
\begin{equation*}
\operatorname{dim}_{\mathcal{H}}D_{\alpha} = \left\{\begin{array}{ll}
4\alpha - 2, & \qquad\alpha \in [1/2,3/4], \\
0, & \qquad \alpha = 3/2, \\
-\infty, & \qquad \text{otherwise},
\end{array} \right.
\end{equation*}
and he gave a characterisation of the sets $D_{\alpha}$ when $1/2 < \alpha < 3/4$ in terms of the rate of convergence of continued fraction approximations. Indeed, these sets are exclusively formed by points $t_{\rho}$ where $\rho$ is irrational.  Let $\{p_n/q_n\}_{n\in\mathbb{N}}$ be the sequence of continued fraction approximations to $\rho$ and define $\gamma_n > 2$ and $\gamma(\rho)$ as
\[ \left| \rho - \frac{p_n}{q_n} \right| = \frac{1}{q_n^{\gamma_n}}, \qquad \gamma(\rho) = \limsup_{n \to \infty}\{\gamma_n \mid q_n \equiv 0,1,3 \pmod{4} \}. \]
Then, $\alpha(t_{\rho}) = 1/2 + 1/(2\gamma(\rho))$. 
By Lemma~\ref{LemmaBoundForContinuedFractions}, this relationship 
allows to reduce the radii of the balls in \eqref{CoverForIrrationals} to cover $\phi(D_{\alpha})$ and therefore to prove the following result.
\begin{thm}\label{TheoremDimensionGeneralised}
Let $\phi$ be Riemann's non-differentiable function in \eqref{Phi}. Then,
\[ \operatorname{dim}_{\mathcal{H}}\phi(D_{\alpha}) \leq \operatorname{dim}_{\mathcal{H}}\bigcup_{\beta \leq \alpha}\phi(D_{\beta}) \leq \frac{4\alpha - 2}{\alpha}  \]
for every $\alpha \in [1/2,3/4]$.
\end{thm}

\section{Tangents}\label{S_Tangents}

Since $\phi(\mathbb{R})$ is expected to be a representative of temporal trajectories in a binormal flow experiment, one may wonder if it has well-defined directions in some sense. Studying whether it has tangents or not from the geometric perspective seems most natural. The main result in this section is the following.

\begin{thm}\label{TheoremTangents}
Riemann's non differentiable function $\phi$, interpreted as a geometric trajectory represented by its image $\phi(\mathbb{R})$, has nowhere a tangent.
\end{thm}

This result is remarkable in two aspects. Hardy \cite{Hardy} and Gerver \cite{Gerver,Gerver2} showed that $\phi$ is differentiable only in $t_{p/q}$ where $p/q$ is an irreducible rational with $q \equiv 2 \pmod{4}$, so one expects to have a tangent in some sense. However, a spiral-like pattern can be seen in the point $\phi(t_{1/2})$ in the right-hand side of Figure~\ref{FigurePhi}. The same happens in all $\phi(t_{p/q})$ alike. This phenomenon arises because $\phi'(t_{p/q})=0$, a cancellation very particular of $\phi$ which does not happen in the case of $\phi_D$. 
For the rest of rationals, $\phi$ is not differentiable in $t_{p/q}$ and indeed no tangents exist in $\phi(t_{p/q})$. However, these points seem to have well-defined right-sided and left-sided tangents, which differ by a right angle. This is suggested by the right-angled corners clearly appreciated in Figure~\ref{FigurePhi}. Last, no immediate conclusion for irrationals can be deduced from Figure~\ref{FigurePhi}.

The concept of tangent has to be determined when working with irregular sets like $\phi(\mathbb{R})$, the image of an almost nowhere differentiable function. In the case of such \textit{irregular} sets $F \subset \mathbb{R}^n$, geometric measure theory supplies useful definitions. In plain terms, if $\operatorname{dim}_{\mathcal{H}}F = s$, $F$ is said to have a tangent in $x \in F$ in direction $\mathbb{V} \subset \mathbb{S}^{n-1}$ if it is concentrated in a cone centred at $x$ with direction $\mathbb{V}$ and arbitrarily small opening angle, when being close enough to $x$. Concentration is measured in terms of the $s$-Hausdorff measure. Since the Hausdorff dimension of $\phi(\mathbb{R})$ is unknown, we adapt this definition using the $1$-Hausdorff content.

\begin{figure}[h]
\centering
\includegraphics[width=0.4\linewidth]{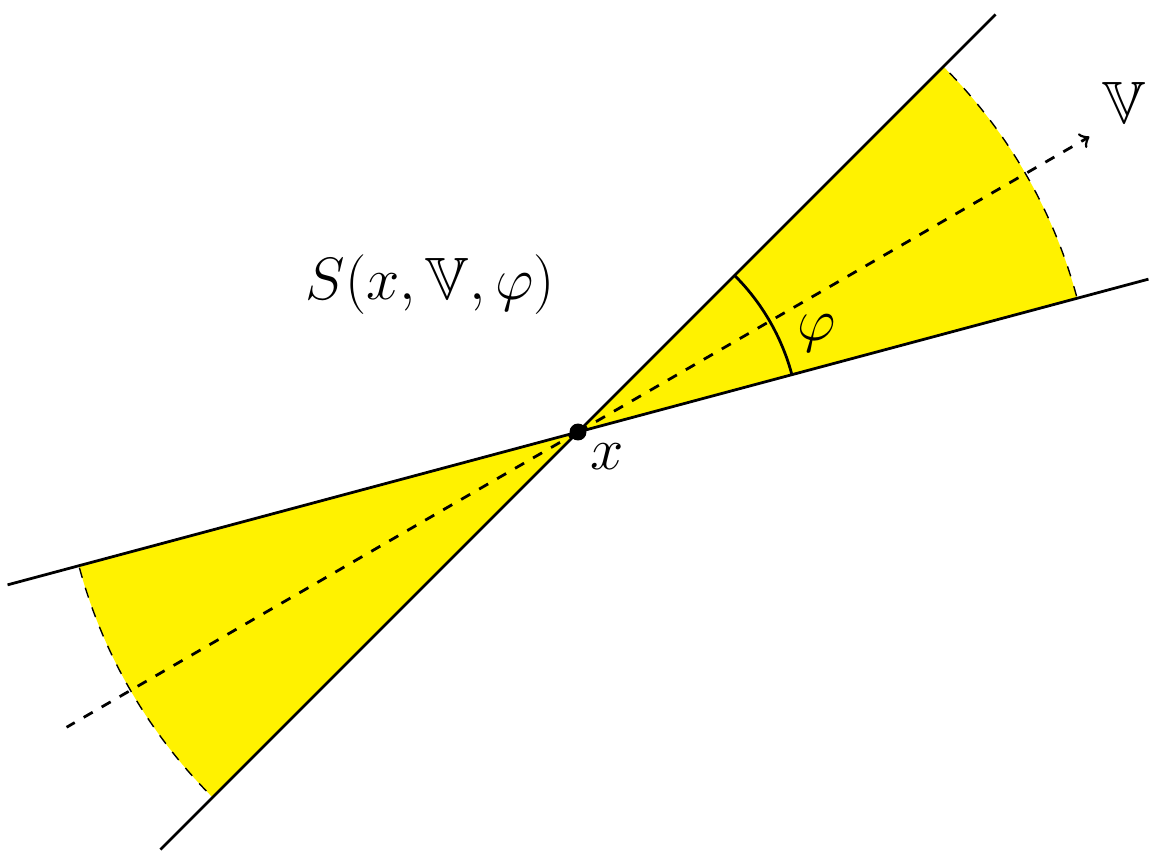}
\caption{The cone $S(x, \mathbb{V},\varphi)$. }
\label{FigureCone1}
\end{figure}

\begin{defin}\label{DefinitionTangents1}
The 1-Hausdorff content of a set $F \subset \mathbb{R}^2$ is 
\[ \mathcal{H}^1_{\infty}(F) = \inf\left\{  \sum_{i}\operatorname{diam}U_i \,  : \, F \subset \bigcup_i U_i, \quad \{U_i\}_i \text{ a countable cover of } F    \right\}.  \]
Let $x \in \mathbb{R}^2$, $\mathbb{V} \in \mathbb{S}^1$ and $\varphi >0$, and denote by $S(x,\mathbb{V},\varphi)$ the closed double cone with vertex in $x$, direction $\mathbb{V}$ and opening angle $\varphi$ (see Figure~\ref{FigureCone1}). Let $t \in \mathbb{R}$. We say that $\mathbb{V}$ is a tangent of $\phi(\mathbb{R})$ in $\phi(t)$ if
\[  \forall \varphi >0, \qquad \lim_{h \to 0}{ \frac{\mathcal{H}^1_{\infty}\left(  \left( \phi(\mathbb{R}) \cap B(\phi(t),h) \right) \setminus S(\phi(t),\mathbb{V},\varphi)  \right)}{h} } = 0.  \]
\end{defin}

However, to prove Theorem~\ref{TheoremTangents} we use an auxiliary, more convenient definition based in the parametrisation $\phi$. Moreover, this approach grasps the right-sided and left-sided tangent phenomenon.

\begin{defin}\label{DefinitionTangent2}
We say that $\mathbb{V} \in \mathbb{S}^1$ is a tangent of $\phi$ in $t\in\mathbb{R}$ on the right if
\[ \lim_{n \to \infty} \frac{\phi(t+h_n) - \phi(t)}{\left| \phi(t+h_n) - \phi(t) \right|} = \mathbb{V} \]
 for every sequence $h_n >0$ converging to zero for which $\phi(t+h_n) - \phi(t) \neq 0$. It is a tangent on the left if the same holds for sequences $h_n <0$. We say $\mathbb{V}$ is a tangent of $\phi$ in $t$ is it is a tangent both from the right and from the left.
\end{defin}

Theorem~\ref{TheoremTangents} is therefore a consequence of the following two results. We remark that Proposition~\ref{PropositionTangents} confirms the guess based in Figure~\ref{FigurePhi} done in the beginning of the section.
\begin{lem}\label{LemmaTangents}
Let $\mathbb{V} \in \mathbb{S}^1$. If $\mathbb{V}$ is not a tangent of $\phi$ in $t \in \mathbb{R}$ as in Definition~\ref{DefinitionTangent2}, then it is not a tangent of $\phi(\mathbb{R})$ in $\phi(t)$ as in Definition~\ref{DefinitionTangents1}.  
\end{lem}
\begin{prop}\label{PropositionTangents}
Let $p/q$ be an irreducible rational number. 
\begin{itemize}
	\item If $q \equiv 0,1,3 \pmod{4}$, there exists an eighth root of unity $e_{p/q} \in \mathbb{C}$ such that 
	\[ \lim_{h \to 0^+}\frac{\phi(t_{p/q} + h) - \phi(t_{p/q})}{\left| \phi(t_{p/q} + h) - \phi(t_{p/q}) \right|} = e_{\frac{p}{q}}\,\frac{1+i}{\sqrt{2}}, \qquad \lim_{h \to 0^+}\frac{\phi(t_{p/q} - h) - \phi(t_{p/q})}{\left| \phi(t_{p/q} - h) - \phi(t_{p/q}) \right|} = e_{\frac{p}{q}}\,\frac{1-i}{\sqrt{2}}  \]
			
	\item If $q \equiv 2 \pmod{4}$, then for every $\mathbb{V} \in \mathbb{S}^1$ there exist sequences $h_n,r_n \to 0^+$ such that 
	\[ \lim_{n \to \infty}\frac{\phi(t_{p/q} + h_n) - \phi(t_{p/q})}{\left| \phi(t_{p/q} + h_n) - \phi(t_{p/q}) \right|} = \mathbb{V} = \lim_{n\to\infty}\frac{\phi(t_{p/q} - r_n) - \phi(t_{p/q})}{\left| \phi(t_{p/q} - r_n) - \phi(t_{p/q}) \right|}. \]
\end{itemize}
Let $\rho$ be an irrational number. Then, there exists an open set $V \subset \mathbb{S}^1$ such that for every direction $\mathbb{V} \in V$, there is a sequence $h_n \to 0$ such that 
\[   \lim_{n \to \infty}\frac{\phi(t_{p/q} + h_n) - \phi(t_{p/q})}{\left| \phi(t_{p/q} + h_n) - \phi(t_{p/q}) \right|} = \mathbb{V}.  \]
\end{prop}

The proof of Lemma~\ref{LemmaTangents} is based on the observation that the 1-Hausdorff content $\mathcal{H}^1_{\infty}$ of a piece of a curve is essentially its diameter. On the other hand, Proposition~\ref{PropositionTangents}, like Lemma~\ref{LemmaBoundForContinuedFractions}, is a consequence of the asymptotic behaviour of $\phi$ at points $t_{p/q}$ corresponding to rational points. Duistermaat \cite{Duistermaat} computed (in terms of $\phi_D$) the asymptotic behaviour of $\phi( t_{p/q}+h)-\phi(t_{p/q})$ up to the error term $h^{3/2}$, which is enough for rational points with $q\equiv 0,1,3 \pmod{4}$. For the rest of rational and irrational numbers, the precise expression of the $h^{3/2}$ term is needed
and is computed in \cite{Eceizabarrena1} following the ideas of \cite{Duistermaat}. 
The case of rational numbers with $q\equiv 2 \pmod{4}$ follows immediately. Finally, as shown in \cite{Eceizabarrena2}, the irrational points are treated by their approximations by means of continued fractions.

\section*{Acknowledgements}
The author is thankful to Valeria Banica, Albert Mas, Xavier Tolsa and Luis Vega for their suggestions and help. 

This work was supported by Spain's Ministry of Education, Culture and Sport  [FPU15/03078], the ERCEA [Advanced Grant 2014 669689 - HADE], the Basque Government [BERC 2018-2021] and the Ministry of Science, Innovation and Universities [BCAM Severo Ochoa accreditation SEV-2017-0718].


\begin{thebibliography}{00}


\bibitem{ChamizoCordoba}
F. Chamizo, A. C\'ordoba, Differentiability and dimension of some fractal Fourier series. Adv. Math. 142 (1999), no.2, 335-354. \href{https://doi.org/10.1006/aima.1998.1792}{https://doi.org/10.1006/aima.1998.1792}

\bibitem{ChamizoUbis2007}
F. Chamizo, A. Ubis, Some Fourier series with gaps. J. Anal. Math. 101 (2007) 179-197. \href{https://doi.org/10.1007/s11854-007-0007-z}{https://doi.org/10.1007/s11854-007-0007-z}

\bibitem{ChamizoUbis2014}
F. Chamizo, A. Ubis, Multifractal behaviour of polynomial Fourier series. Adv. Math. 250 (2014) 1-34. \href{https://doi.org/10.1016/j.aim.2013.09.015}{https://doi.org/10.1016/j.aim.2013.09.015}

\bibitem{DeLaHozVega}
F. De la Hoz, L. Vega, Vortex Filament Equation for a Regular Polygon. Nonlinearity 27 (2014) 3031-3057. \href{https://doi.org/10.1088/0951-7715/27/12/3031}{https://doi.org/10.1088/0951-7715/27/12/3031}

\bibitem{Duistermaat}
J.J. Duistermaat, Selfsimilarity of Riemann's nondifferentiable function. Nieuw Arch. Wisk. (4) 9 (1991) no.3, 303-337.

\bibitem{Eceizabarrena1}
D. Eceizabarrena. Asymptotic behaviour and Hausdorff dimension of Riemann's non-differentiable function. Preprint, \href{https://arxiv.org/abs/1910.02530}{arXiv:1910.02530}.

\bibitem{Eceizabarrena2}
D. Eceizabarrena. Geometric differentiability of Riemann's non-differentiable function. Preprint, \href{https://arxiv.org/abs/1910.02536}{arXiv:1910.02536}.

\bibitem{Gerver}
J. Gerver, The differentiability of the Riemann function at certain rational multiples of $\pi$. Amer. J. Math. 92 (1970) 33-55. \href{https://doi.org/10.2307/2373496}{https://doi.org/10.2307/2373496}

\bibitem{Gerver2}
J. Gerver, More on the differentiability of the Riemann function. Amer. J. Math. 93 (1971) 33-41. \href{https://doi.org/10.2307/2373445}{https://doi.org/10.2307/2373445}

\bibitem{Hardy}
G.H. Hardy, Weierstrass' non-differentiable function. Trans. Amer. Math. Soc. 17 (1915), no.3,  301-325. \href{https://doi.org/10.2307/1989005}{https://doi.org/10.2307/1989005}

\bibitem{HardyLittlewood}
G.H. Hardy, J.E. Littlewood, Some problems of Diophantine approximations (II). Acta Math. 37 (1914) 193-239. \href{https://doi.org/10.1007/BF02401834}{https://doi.org/10.1007/BF02401834}

\bibitem{HolschneiderTchamitchian}
M. Holschneider, P. Tchamitchian, Pointwise analysis of Riemann's "nondifferentiable" function. Invent. Math. 105 (1991), no.1, 157-175. \href{https://doi.org/10.1007/BF01232261}{https://doi.org/10.1007/BF01232261}

\bibitem{Jaffard}
S. Jaffard, The spectrum of singularities of Riemann's function. Rev. Mat. Iberoamericana, 12 (1996), no.2, 441-460. \href{https://doi.org/10.4171/RMI/203}{https://doi.org/10.4171/RMI/203}


\bibitem{Pastor}
C. Pastor, On the regularity of fractional integral of modular forms. Trans. Amer. Math. Soc. 372 (2019), no.2, 829-857. \href{https://doi.org/10.1090/tran/7418}{https://doi.org/10.1090/tran/7418}

\bibitem{Weierstrass}
K. Weierstrass, \"Uber continuirliche functionen eines reellen arguments, die f\"ur keinen werth des letzteren einen bestimmten differentialquotienten besitzen. Mathematische Werke II, K\"onigl. Akad. Wiss. (1895) 71-74.


\end{thebibliography}
\end{document}